\documentclass[12pt,draftcls,onecolumn]{IEEEtran}

\IEEEoverridecommandlockouts                              

\usepackage{graphics} 
\usepackage{epsfig} 
\usepackage{mathptmx} 
\usepackage{times} 
\usepackage{amsmath} 
\usepackage{amssymb}  
\usepackage{soul}
\usepackage[mathcal]{euscript} 
\usepackage{mathtools}
\usepackage{mathrsfs}

\usepackage{graphicx}
\graphicspath{{Figures/}} 
\usepackage{caption} 
\usepackage{subfig} 

\usepackage{url}

\usepackage{xcolor}
\definecolor{USred}{rgb}{0.74,0.1,0.1}
\definecolor{USblue}{rgb}{0.2,0.2,0.7}
\definecolor{green1}{cmyk}{0.82,0,1,0.3}

\definecolor{Royalblue}{cmyk}{1,0.30,0.2,0.2}

\newcommand{\numberset}{\mathbb}
\newcommand{\NN}{\numberset{N}}

\newcommand{\RR}{\numberset{R}}

\newcommand{\ZZ}{\numberset{Z}}
\newcommand{\EE}{\mathbb{E}}

\newcommand{\B}{\mathcal{B}}
\newcommand{\C}{\mathcal{C}}

\newcommand{\G}{\mathcal{G}}

\newcommand{\N}{\mathcal{N}}
\newcommand{\Q}{\mathcal{Q}}

\newcommand{\argmin}{\operatornamewithlimits{argmin}}
\newcommand{\argmax}{\operatornamewithlimits{argmax}}
\renewcommand{\vec}{\boldsymbol}
\DeclareMathOperator{\tr}{tr}

\DeclareRobustCommand{\vect}[1]{
	\ifcat#1\relax
	\boldsymbol{#1}
	\else
	\mathbf{#1}
	\fi}

\newcommand\indep{\protect\mathpalette{\protect\independenT}{\perp}}
\def\independenT#1#2{\mathrel{\rlap{$#1#2$}\mkern2mu{#1#2}}}

\newtheorem{definition}{Definition}
\newtheorem{theorem}{Theorem}
\newtheorem{proposition}{Proposition}
\newtheorem{problem}{Problem}
\newtheorem{remark}{Remark}
\newtheorem{example}{Example}

\title{Identification of Sparse Reciprocal Graphical Models}

\author{Daniele Alpago, Mattia Zorzi, Augusto Ferrante 
\thanks{}
\thanks{D. Alpago, M. Zorzi and A. Ferrante are with the Department of Information Engineering, University of Padova, Padova, Italy; email:	 
	    {\tt\small alpagodani@dei.unipd.it} (D. Alpago)
        {\tt\small zorzimat@dei.unipd.it} (M. Zorzi)
        {\tt\small augusto@dei.unipd.it} (A. Ferrante)}%
\thanks{}%
}

\begin{document}

\maketitle
\thispagestyle{empty}
\pagestyle{empty}

\begin{abstract}
In this paper we propose an identification procedure of a sparse graphical model associated to a Gaussian stationary stochastic process. The identification paradigm exploits the approximation of autoregressive processes through reciprocal processes in order to improve the robustness of the identification algorithm, especially when the order of the autoregressive process becomes large. We  show that the proposed paradigm leads to a regularized, circulant matrix completion problem whose solution only requires computations of the eigenvalues of matrices of dimension equal to the dimension of the process.
\end{abstract}


\section{INTRODUCTION}\label{sec:intro}

\IEEEPARstart{I}{n} the last decade miniaturisation led to an ubiquitous pervasiveness of technology. As a consequence, the number of available high-dimensional data is skyrocketing in every scientific and applicative domain. Several methods to deal with problems involving high-dimensional data have been recently proposed in the literature \cite{BDS}. In this paper we are focusing on graphical models, that represent a possible tool to deal with high-dimensionality of the data \cite{Lauritzen}. Graphical representations provide an immediate visual intuition on the data interdependence. The simplest  graphical model is an undirected graph that can be associated with a Gaussian random vector
\cite{candes2012exact,candes2010matrix,ZorzCh}: nodes correspond to the components of the random vector, and there is an edge between two nodes if the corresponding components are conditionally dependent given all the others.
Very often data are given as time-series and can thus be modelled as stochastic processes. Also such processes  can  be represented as graphical models: in \cite{songsiri2010graphical} a maximum likelihood approach has been proposed for graphical model estimation of autoregressive (AR) Gaussian processes, exploiting the fact that conditional independence relations translates in zero entries in the inverse of the power spectral density. Since the conditional independent pairs are not known a priori, a sparsity-inducing regularizer can be introduced in the problem \cite{SongVan}, leading to a graphical model with a \emph{sparse} structure, i.e. with few edges. The reduced number of edges gives a double advantage: the graph gives a clearer representation of the phenomenon we are observing, and the number of parameters that have to be estimated is reduced. Moreover, sparsity makes the identification procedure less subject to overfitting, thus leading to a parsimonious estimated model. In \cite{LindqARMA} the aforementioned paradigm has been extended to ARMA Gaussian processes wherein the moving average (MA) part is introduced by a scalar ``prior'' power spectral density. The MA part can be determined from logarithmic moments (cepstral coefficients) of the spectrum \cite{PICCI_LINDQUIST, BEL01, BEL02, MK85} but it is not clear how to impose such moments together with the constraints on the covariance lags, so that the estimated spectrum reflects the underlying graphical structure.
The paradigm in \cite{SongVan} can be extended to latent-variable graphical models \cite{ChaParWill,ZorzSep}.

The identification of Gaussian ARMA graphical models can be
performed, in principle, by using the method in \cite{SongVan}: indeed, an ARMA process is well approximated by an high-order AR process. The problem, however, is that the optimization procedure involves the inversion and the eigenvalue decomposition of matrices whose dimensions are proportional to the product of the order of the 
AR process by the dimension of the data. As a consequence, the procedure becomes numerically less robust when the process is high-dimensional and the AR approximation is sufficiently accurate and hence ``long''.

In this paper we consider the problem of identifying sparse graphical models for Gaussian reciprocal processes defined in the ``discrete circle'' ${\ZZ}_N$ (the group of the integers modulo $N$). The latter constitute a particular class of periodic processes, \cite{CarFerrPav}, \cite{Rec1}, \cite{Rec2}, \cite{Rec3}, \cite{ringh2016multidimensional}, \cite{PICCI_LINDQUIST}.  It is possible to show that a Gaussian AR process of order $n$ can be well approximated by a reciprocal process of the same order $n$ as long as the  period $N$ is sufficiently large \cite{FerrLinAlg}, \cite{PICCI_LINDQUIST}. 
We will show that the identification problem involves block-circulant matrices. This is a big numerical advantage because the inversion and the eigenvalue decomposition of such matrices can be performed robustly even if $n$ is large \cite{ringh2015fast}. Accordingly, the proposed  paradigm can represent potentially a robust method for estimating ARMA graphical models.


The paper is organized as follows: In Section \ref{sec:recpr} we recall the fundamental results regarding the identification of reciprocal  processes. In Section \ref{sec:gms} we introduce graphical models for the reciprocal processes. In Section \ref{sec:id} we present the problem of estimating a sparse graphical model for a reciprocal process and in Section \ref{sec:ne} we show how the proposed method behaves in a numerical example. Finally, in Section \ref{sec:concl} we draw the conclusions.

\emph{Notation and background.} We denote by $\NN,\,\ZZ,\,\RR,\,\RR^{p\times q}$ the set of natural, integers, real numbers and $p\times q$ real matrices, respectively. For any sub-interval $(x_1,x_2):=\{x:\,x_1<x<x_2\}$ of an interval $(a,b)\subset\RR$, we denote with $(x_1,x_2)^c$ the complement set of $(x_1,x_2)$ in $(a,b)$. Given a matrix $A\in\RR^{p\times q}$, we denote by $A^\top$ its transpose and by $\ker(A)$ its kernel, while $I_p$ denotes the identity matrix of order $p$. If $A\in\RR^{p\times p}$ is a square matrix, $\tr(A),\,\det(A)$ and $A^{-1}$ denote, respectively, the trace of $A$, the determinant of $A$ and its inverse, while $\text{diag}(A)\in\RR^p$ is the vector whose entries are the diagonal elements of $A$. If $A$ is symmetric, $A>0$ and $A\ge 0$ indicate that it is positive definite or positive semi-definite, respectively. Moreover, we will use $\EE[\cdot]$ to denote the expectation operator.

In this paper $n$ will always represent the order of the original AR process which will be of dimension $m$. We will denote by $N$  the length the period of the reciprocal 
process. Thus $n$, $m$ and $N$ are fixed numbers.  For simplicity we  assume that $N\in\NN$ is an even number, however the results in this paper can be easily adapted to the case when $N$ is odd.
Also, we assume $N>2n$.
A central role in this paper will be played by the vector space $\C\subset\RR^{mN\times mN}$ containing all the (real) symmetric, block-circulant matrices
\[
   \vect{C} = \text{circ}\{C_0,C_1,\dots,C_{\frac{N}{2}-1},C_{\frac{N}{2}},C_{\frac{N}{2}-1}^\top,\dots,C_1^\top\},
\]
whose first block-column is composed by the $m\times m$ blocks $C_0,C_1,\dots,C_{\frac{N}{2}-1},C_\frac{N}{2},C_{\frac{N}{2}-1}^\top,\dots,C_1^\top$. For any $\vect{C},\,\vect{D}\in\C$, the inner product on $\C$ is defined by $\left<\vect{C},\vect{D}\right>_\C:= \tr(\vect{C}^\top\vect{D})$. We define the \emph{symbol} of the circulant matrix $\vect{C}\in\C$ as the $m\times m$ pseudo-polynomial
\begin{equation}\label{eq:symb}
\Phi(\zeta) := \sum_{k=0}^{N-1}\,C_k\,\zeta^{-k},\quad \text{ with } \quad C_k = C_{N-k}^\top\text{ for } k>\frac{N}{2},
\end{equation}\noindent
where $\zeta:=e^{i\frac{2\pi}{N}}$ is the $N$-th root of unity. It is useful to recall the following result on block-circulant matrices.
\begin{proposition}\label{prop:diagsymb}
	Let $\vect{C}$ be a block-circulant matrix with symbol $C(\zeta)$ defined by \eqref{eq:symb}. Then
	\begin{equation}
	\vect{C} =
	\vect{F}^*
	\text{diag}\left\{\Phi(\zeta^0),\,\Phi(\zeta^1),\,\dots,\,\Phi(\zeta^{N-1})\right\}
	\vect{F},
	\end{equation}
	where $\vect{F}$ is the $mN\times mN$ (Fourier) unitary matrix
	{\small\begin{equation*}
	\vect{F}=
	\frac{1}{\sqrt{N}}
	\begin{bmatrix}
	\zeta^{-0 \cdot 0}I_m    & \zeta^{-0 \cdot 1}I_m     & \ldots & \zeta^{-0 \cdot (N-1)}I_m    \\
	\zeta^{-1 \cdot 0}I_m     & \zeta^{-1 \cdot 1}I_m     & \ldots & \zeta^{-1 \cdot (N-1)}I_m     \\
	\vdots                   & \vdots                   & \ddots & \vdots                       \\
	\zeta^{-(N-1) \cdot 0}I_m & \zeta^{-(N-1) \cdot 1}I_m & \ldots & \zeta^{-(N-1) \cdot (N-1)}I_m \\
	\end{bmatrix}.
	\end{equation*}}
\end{proposition}
This is a classical result in the scalar case; the general analysis for block-circulant matrices can be found, for instance, in \cite[page 6]{tesi}. We introduce the subspace $\B\subseteq\C$ of symmetric, banded block-circulant $mN\times mN$ matrices of bandwidth $n$, with $N>2n$, containing matrices of the form
\begin{equation}\label{eq:bandmtx}
\vect{B} = \text{circ}\{B_0,B_1,\dots,B_n,0,\dots,0,B_n^\top,\dots,B_1^\top\},
\end{equation}
that inherits the inner product defined on $\C$. Note that, according to definition \eqref{eq:symb}, the symbol of a banded matrix $\vect{B}\in\B$ is
\begin{equation*}
\Psi(\zeta) = \sum_{k=-n}^n\, B_k\,\zeta^{-k}, \qquad B_{-k} = B_k^\top.
\end{equation*} 
If $\mathscr{B} := [0,n]$ denotes the set of indexes of the blocks in the banded structure, then the projection operator $\mathsf{P}_\mathscr{B}:\C\to\B$ is defined as 
\[
\mathsf{P}_\mathscr{B}(\vect{C}) := \text{circ}\{C_0,C_1,\dots,C_n,0,\dots,0,C_n^\top,\dots,C_1^\top\}.
\]

\section{RECIPROCAL PROCESS IDENTIFICATION}\label{sec:recpr}
Let $\{\vect{y}(k),\,k=1,2,\dots,N\}$, be a zero-mean, $m$-dimensional Gaussian stationary stochastic process defined on a finite interval $[1,N]$. More explicitly we have $\vect{y}(k) := [\vec{y}_1(k)\,\dots\,\vec{y}_m(k)]^\top\in\RR^m$, $k=1,\dots,N$, therefore the process is completely characterized by the random vector $\vect{y} := [\vec{y}_1(1)\,\dots\,\vec{y}_m(1)\,\dots\,\dots\,\vec{y}_1(N)\,\dots\,\vec{y}_m(N)]^\top\in\RR^{mN}$. In \cite{CarFerrPav} it has been shown that $\vect{y}$ is a restriction of a wide-sense stationary periodic process of period $N$ defined on the whole integer line $\ZZ$ if and only if the $mN\times mN$ covariance matrix $\vect{\Sigma}$ of $\vect{y}$ is symmetric block-circulant:
\begin{equation}\label{eq:circmtx}
   \vect{\Sigma} =  \text{circ}\{\Sigma_0,\Sigma_1,\dots,\Sigma_\frac{N}{2},\dots,\Sigma_1^\top\},
\end{equation}
where $\EE[\vect{y}(i)\vect{y}(j)^\top] = \Sigma_{i-j}$, $i,j = 1,\dots,N$, are the covariance lags of the process such that $\Sigma_k = \Sigma_{N-k}^\top$ for $k>N/2$. In view of the above equivalence, we will denote with $\vect{y}$ both the wide-sense stationary periodic process defined in the whole line $\ZZ$ and its restriction, depending on the context. 
A particular class of stationary periodic processes is represented by reciprocal processes.

\begin{definition}
	We say that $\vect{y}$ is a (periodic) reciprocal process of order $n$ on $[1,N]$ if the random variables of the process in the interval $(t_1,t_2)\subset[1,N]$ are conditionally orthogonal to the random variables in $(t_1,t_2)^c$, given the $2n$ boundary values $\vect{y}(t_1-n+1),\dots,\vect{y}(t_1),\vect{y}(t_2),\dots,\vect{y}(t_2+n-1)$, where the sums $t-k$ and $t+k$ are to be understood modulo $N$.
\end{definition}

The following result has been proved in \cite[Theorem 3.3]{CarFerrPav}, and it simply states that a reciprocal model is completely specified by a block-circulant matrix of type \eqref{eq:circmtx} whose inverse is banded, block-circulant as in \eqref{eq:bandmtx}.

\begin{theorem}\label{thm:conc}
	A non-singular $mN \times mN$-dimensional matrix $\vect{\Sigma}$ is the covariance matrix of a periodic reciprocal process of order $n$ if and only if its inverse is a positive definite symmetric block-circulant matrix which is banded of bandwidth $n$, namely $\vect{\Sigma}^{-1}\in\B$.
\end{theorem}
\mbox{\\}

We are now ready to deal with the identification problem of a reciprocal process. Let $\hat{\Sigma}_0,\dots,\hat{\Sigma}_n$ be given estimates of the first $n+1$ covariance lags $\Sigma_0,\dots,\Sigma_n$ of the underlying reciprocal process. According to Theorem \ref{thm:conc}, the identification of a reciprocal model can be formulated as a matrix completion problem.

\begin{problem}\label{pb:covext}
	Given the $n+1$ estimates $\hat{\Sigma}_0,\dots,\hat{\Sigma}_n$, compute a sequence $\Sigma_{n+1},\dots,\Sigma_\frac{N}{2}$, in such a way to form a symmetric, positive definite block-circulant matrix 
	\[
	   \vect{\Sigma} = \text{circ}\{\hat{\Sigma}_0,\dots,\hat{\Sigma}_n,\Sigma_{n+1},\dots,\Sigma_\frac{N}{2},\dots,\Sigma_{n+1}^\top,\hat{\Sigma}_n^\top,\dots,\hat{\Sigma}_1^\top\},
	\]
	with $\Sigma^{-1}\in\B$.
\end{problem}
It has been shown in \cite{CarFerrPav} that the condition $\Sigma^{-1}\in\B$ is equivalent to maximizing the entropy of the process so that the previous problem is equivalent to the following optimization program:
\begin{equation}\label{op:mep}
\begin{aligned}
\argmax_{\vect{\Sigma}\in\C} &\quad \log\det\vect{\Sigma}\\
\text{subject to } &\quad \vect{\Sigma}>0\\
                   &\quad \mathsf{P}_{\mathscr{B}}(\vect{\Sigma}-\hat{\vect{\Sigma}}) = 0.
\end{aligned}
\end{equation}
whose dual problem has been proven to be
\begin{equation}\label{op:mepdu}
\begin{aligned}
\argmin_{\vect{X}\in\B} &\quad -\log\det\vect{X} + \left<\vect{X},\,\hat{\vect{\Sigma}}\right>_\C\\
\text{subject to } &\quad \vect{X}>0
\end{aligned}
\end{equation}
where $\hat{\vect{\Sigma}}\in\B$ is the symmetric, banded block-circulant matrix of bandwidth $n$,
$
   \hat \Sigma= \text{circ}\{\hat \Sigma_0,\hat \Sigma_1,\dots,\hat \Sigma_n,0,\dots,0,\hat \Sigma_n^\top,\dots,\hat \Sigma_1^\top\},
$
containing the covariance lags estimated from the data and
the optimal value of dual variable $\vect{X}$ is indeed equal to $\Sigma^{-1}$, i.e. the inverse of the solution of problem \ref{pb:covext}.
It can be shown that strong duality holds between \eqref{op:mep} and \eqref{op:mepdu}, so that \eqref{op:mep} and \eqref{op:mepdu} are equivalent. In what follows we assume that $\hat{\vect{\Sigma}}>0$ as it is a necessary condition for problem \eqref{op:mep} to be feasible. In the case that $\hat{\vect{\Sigma}}$ is not positive definite, we can consider a positive definite banded block-circulant matrix sufficiently close to $\hat{\vect{\Sigma}}$ which can be obtained by solving a structured covariance estimation problem, see \cite{structcov1}, \cite{structcov2}.

\begin{remark}\label{rmk:approx} 
	Recall that, for $N\to\infty$, Toeplitz matrices can be approximated arbitrarily well by circulant matrices, see \cite[Lemma 4.2]{Gray}. Accordingly, for $N\to\infty$, problem \ref{pb:covext} consists in searching a completion that leads to an infinite positive definite block-Toeplitz covariance matrix, i.e. such that the Fourier transform of the resulting extended sequence is a power spectral density. By Theorem 3.1 of \cite{FerrLinAlg}, for $N\to\infty$, problem \eqref{op:mep} is the classical Burg's maximum entropy problem \cite{BurgPhd,family_zorzi,ZORZI201587,PICCI_LINDQUIST}, whose solution is an autoregressive process of order $n$. In light of this observation, we can understand the reciprocal process solution of \eqref{op:mepdu} as an approximation of the AR process solution of the Burg's maximum entropy problem. In the following sections we will exploit this approximation for the identification of sparse AR graphical models.
\end{remark}

\begin{figure}[h!]\centering
	\includegraphics[scale=0.8]{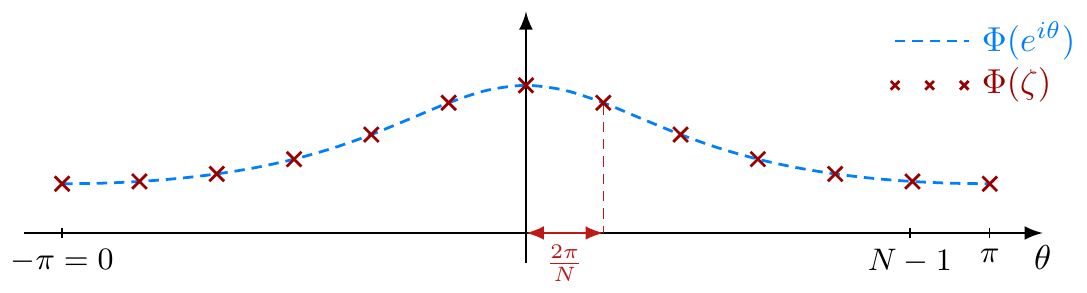}
	\caption{Power spectrum $\Phi(e^{i\theta})$ and its sampled version $\Phi(\zeta)$ with $N=12$ samples.}\label{fig:recapx}
\end{figure}

The reciprocal approximation described in Remark \ref{rmk:approx} has also an interesting interpretation in the frequency domain. Let $\Phi$ denote the power spectrum of the autoregressive, wide-sense stationary process $\vect{x}$,
\begin{equation*}
\Phi(e^{i\theta}) = \sum_{k=-\infty}^{\infty}\,R_k\,e^{-i\theta k},\qquad R_{-k}=R_k^\top,\ \ \theta\in[-\pi,\pi].
\end{equation*}
The reciprocal approximation translates in sampling the power spectrum over the interval $[-\pi,\pi]$, with sample frequency $2\pi /N$, obtaining the \emph{symbol} of the covariance matrix of the corresponding reciprocal process:
\[
\Phi(\zeta) = \sum_{k=0}^{N-1}\, \Sigma_k \,\zeta^{-k},\qquad \Sigma_{k}=\Sigma_{N-k}^\top\text{ for } k>\frac{N}{2}.
\]
Figure \ref{fig:recapx} gives an intuitive idea of this approximation. According to Proposition \ref{prop:diagsymb}, the covariance matrix $\vect{\Sigma}$ of the reciprocal process $\vect{y}$ that approximates $\vect{x}$ writes as
\begin{equation}\label{eq:covrec}
\vect{\Sigma} = \vect{F}^*\text{circ}\{ \Phi(\zeta^0),\,\Phi(\zeta^1),\dots,\Phi(\zeta^{N-1}) \}\vect{F},
\end{equation}
hence, its inverse
\begin{equation}\label{eq:invcovrec}
\vect{\Sigma}^{-1} = \vect{F}^*\text{circ}\{ \Phi(\zeta^0)^{-1},\,\Phi(\zeta^1)^{-1},\dots,\Phi(\zeta^{N-1})^{-1} \}\vect{F},
\end{equation}
can be robustly computed by inverting the $N$ blocks $\Phi(\zeta^0),\,\Phi(\zeta^1),\dots,\Phi(\zeta^{N-1})$, all of size  of size $m\times m$. As a final remark, we recall that eigevalues and eigenvectors of circulant matrices can be robustly computed as well, thanks to the availability of closed-form formulas, see for instance \cite{Gray}.%
\vspace{2mm}
\\
As highlighted by the frequency-domain interpretation, the goodness of the approximation strictly depends on the regularity of the spectrum: the less the spectrum is regular, the larger $N$ has to be chosen in order to get a good approximation of the AR process. 

\section{GRAPHICAL MODELS}\label{sec:gms}
Consider a Gaussian random vector $\vect{x}=[\vec{x}_1\,\dots\,\vec{x}_m]^\top$ with covariance matrix $\vect{\Sigma}$ and let $K := \Sigma^{-1}$ be its \emph{concentration matrix}. The notation 
\[
   \vec{x}_i \indep \vec{x}_j \mid \{\vec{x}_k\}_{k\ne i,j}
\]
means that the random variable $\vec{x}_i$ is \emph{conditionally independent} from the random variable $\vec{x}_j$ given the remaining random variables $\vec{x}_k$, $k\ne i,j$. It can be proven that, \cite{Lauritzen}:
\begin{equation}\label{eq:cic}
	\vec{x}_i \indep \vec{x}_j \mid \{\vec{x}_k\}_{k\ne i,j}
	\iff
	k_{ij}=0,
\end{equation}
where $k_{ij}:=(K)_{ij}$ is the element in position $(i,j)$ in the concentration matrix $K$, $i,\,j=1,\dots,m$. Relation \eqref{eq:cic} defines an undirected graph $\G = (V,E)$, $E\subset V\times V$, associated to the random vector $\vect{x}$, whose nodes are the components $\vec{x}_1,\dots,\vec{x}_m$ of $\vect{x}$, and the absence of edges describes conditional independence  between the components, namely for $i\ne j$, 
\[
    (i,j)\notin E \iff \vec{x}_i \indep \vec{x}_j \mid \{\vec{x}_k\}_{k\ne i,j}.
\]
The graph $\G$ is called the \emph{graphical model} associated to $\vect{x}$. Property \eqref{eq:cic} provides a complete characterization of the graphical model associated to a certain Gaussian random vector in terms of its concentration matrix.
In practice, there is a large interest in \emph{sparse graphical models}, i.e. graphs that describe the interactions between a large number of components $x_i$s with few edges (equivalently with $K$ being a sparse matrix), and thus give an easily understandable description of the underlying system we are modeling.

Although there is a large literature that deals with sparse graphical models \cite{SongVan}, \cite{ZorzSep,ZorzCh,songsiri2010graphical,LindqARMA}, \cite{ChaParWill,candes2012exact,candes2010matrix}, the problem of deriving such models for the case in which the underlying process is a reciprocal process seems to not have been considered till now. In what follows we present sparse graphical models associated to reciprocal processes introduced in Section \ref{sec:recpr}. We will describe how the combination of the underlying reciprocal structure and the sparsity constraint on the concentration matrix of the process impact the properties of the resulting sparse graphical model. The proofs of the following results can be found in \cite{tesi}.\\
Let $\vect{y}$ be a Gaussian, periodic, reciprocal process of order $n$ defined on $[1,N]$ with covariance matrix $\vect{\Sigma}\in\C$ and let $\vect{S}:=\vect{\Sigma}^{-1}$ denote its \emph{concentration matrix} so that, according to Theorem \ref{thm:conc},
\begin{equation}\label{eq:Smtx}
	\vect{S} =\text{circ}\{S_0,S_1,\dots,S_n,0,\dots,0,S_n^\top,\dots,S_1^\top\}.
\end{equation}
In the following we will generalize the characterization of conditional independence we have given in the classical setting of Gaussian random vectors to the case of graphical models associated to Gaussian reciprocal processes. For this purpose, it is useful to define the $j$-th component of the process $\vect{y}$ as the $\RR^N$-valued vector $\vect{y}_j:=[\vec{y}_j(1) \dots \vec{y}_j(N)]^\top$, obtained by stacking all the $j$-th components of the process for each $k=1,\dots,N$. The components of the reciprocal process are defined \emph{for any} $k\in\ZZ$. The process, however, is periodic of period $N$ so that we can impose conditional independence only for $k\in[1,N]$. By Property \eqref{eq:cic}, this implies that $\vect{S}$ is a \emph{sparse matrix} and that the blocks $S_0, S_1,\dots,S_n$ have common support $\Omega\subseteq\{(i,j):\,i,j=1,\,\dots,\,m\}$ namely,
\begin{equation}\label{eq:cs}
(S_k)_{ij} = (S_k)_{ji} = 0,\qquad k=0,\,\dots\,,n,\quad\forall\,(i,j)\in\Omega^c,
\end{equation}
where $\Omega$ is the set of conditionally dependent pairs that necessarily contains all the pairs $(i,\,i),\,i=1,\dots,m$, since conditional independence is not defined between one variable and itself. The above relation is equivalent to
\begin{equation}\label{eq:cs1}
\begin{aligned}
\EE\big[\vec{y}_i(t_1)\,\vec{y}_j(t_2)\mid\,&\vec{y}_h(s),\,h\ne i,j,\,\,s=1,\dots,N,\\
& \vec{y}_i(s_1),\,s_1\ne t_1,\,\vec{y}_j(s_2),\,s_2\ne t_2\big]=0,
\end{aligned}
\end{equation}
for any $t_1,\,t_2\in[1,N]$ and for any pair $(i,j)\in\Omega^c$. The following result has been proved in \cite{tesi}.

\begin{proposition}\label{thm:tesi}
	Condition \eqref{eq:cs1} is equivalent to
	\begin{equation}\label{eq:condindp}
	\EE\left[ \vec{y}_i(t_1)\,\vec{y}_j(t_2)\mid 
	\vec{y}_h(s),\,h\ne i,j,\,\,s=1,\dots,N\right] = 0,
	\end{equation}
	for any $t_1,\,t_2\in[1,N]$ and for any $(i,j)\in\Omega^c$.
\end{proposition}

The above result reflects the fact that the random variables $\{\vec{y}_i(s_1),\vec{y}_j(s_2),\,s_1\ne t_1,\,s_2\ne t_2\}$ do not play any role in the conditioning \eqref{eq:cs1}.
According to \eqref{eq:condindp} we can associate a graphical model to the process $\vect{y}$, whose nodes are its $m$, $N$-dimensional components $\vect{y}_1,\dots,\vect{y}_m$ and the presence of an edge between two vectors $\vect{y}_i$ and $\vect{y}_j$, $i\ne j$, means that $\vect{y}_i$ and $\vect{y}_j$ are conditionally dependent. According to the characterization of conditional dependence, such an edge is described by the quantities 
\[
\EE\left[ \vec{y}_i(t_1)\,\vec{y}_j(t_2)\mid 
\vec{y}_h(s),\,h\ne i,j,\,\,s=1,\dots,N\right]
\]
for any $t_1,\,t_2\in[1,N]$ and for $i,j=1,\dots,m$, i.e. the edge $(i,j)$ is completely characterized by the vector
\[
\left[ (S_0)_{ij}\,\,(S_1)_{ij}\,\,\dots\,\,(S_n)_{ij}\,\, 0 \,\, \dots \,\, 0 \,\, (S_n)_{ji}\,\,\dots\,\, (S_1)_{ji}\right].
\]
By Property \eqref{eq:cic}, the identification of such graphical models consists in the estimation of the sequence $S_0,S_1,\dots,S_n$ obeying the structural constraint (\ref{eq:cs}).

\begin{example}
	Consider the case in which $N=6$, $n=2$ and $m=4$, and suppose that the graphical model associated to process $\vect{y}$ is the one depicted in
	Figure \ref{fig:graphprop}.
	\begin{figure}[h!]\centering
		\includegraphics[scale=1.2]{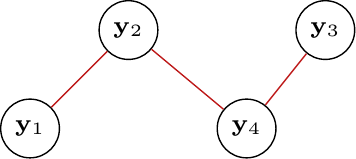}
		\caption{Example of a sparse reciprocal graphical model.}\label{fig:graphprop}
	\end{figure}\\
    In this case, the concentration matrix of vector $\vect{y}$ is
    \[
       \vect{S} =\text{circ}\{S_0,\,S_1,\,S_2,\,0,\,S_2^\top,\,S_1^\top\},
    \]
    where $(S_j)_{12},\,(S_j)_{24},\,(S_j)_{34}$, $j=0,1,2$, are the unique entries different from zero. In this case, $\Omega=\{ (i,i):\,i=1,\dots,4\}\cup\{(1,2),\,(2,4),\,(3,4)\}$.
\end{example}

\section{IDENTIFICATION OF SPARSE RECIPROCAL GRAPHICAL MODELS}\label{sec:id}
Let $\vect{x}:=\{\vect{x}(t):\,t\in\ZZ\}$ be an $m$-dimensional, AR Gaussian stationary process of order $n$, 
\begin{equation}\label{eq:AR}
\sum_{k=0}^{n}\,B_k\,\vect{x}(t-k) = \vect{e}(t),\qquad\vect{e}(t)\sim\N(0,I_m),
\end{equation}
defined in the whole integer line $\ZZ$, and let $R_k := \EE[\vect{x}(t)\vect{x}(t-k)^\top]$, $k\in\ZZ$, be its $k$-th covariance lag. Suppose now that $T$ observations $x(1),\dots,x(T)$ are available, and let
\begin{equation}\label{eq:estcov}
\hat{R}_k = \frac{1}{T}\sum_{t=k}^T\,x(t)x(t-k)^\top,\qquad k = 0,1,\dots,n,
\end{equation}
be estimates of the first $n+1$ covariance lags $R_0,\dots,R_n$.
In view of Remark \ref{rmk:approx}, the idea is to approximate process $\vect{x}$ with a Gaussian reciprocal process $\vect{y}$ of order $n$ defined over the interval $[1,N]$, with $N$ sufficiently large, having a sparse graphical model.

By what we have explained in Section \ref{sec:recpr}, the introduction of the reciprocal approximation allows to obtain a robust procedure even in the case that $n$ is large. In fact, since the matrices that are involved in the optimization are symmetric and block-circulant, according to relations \eqref{eq:covrec} and \eqref{eq:invcovrec}, we can compute the corresponding inverse matrices and eigenvalues in a  robust way. Moreover, it is apparent from \eqref{eq:invcovrec} that the identification algorithm we are proposing scales with respect to $n$ because the dimensions of the matrices, whose eigenvalues must be computed in the optimization procedure, depend only on $m$. This ensures robustness in the results even if the order of the AR process we are considering is large.\\
Now we can formally state the identification problem. 

\begin{problem}\label{prob2}
Consider an $m$-dimensional process $\vect{x}$ and let $\hat{R}_0,\dots,\hat{R}_n$ be the estimates of the first $n+1$ covariance lags of $\vect{x}$ given by \eqref{eq:estcov}. Set $\Sigma_0:=\hat{R}_0,\dots,\Sigma_n:=\hat{R}_n$. Compute the blocks $\Sigma_{n+1},\dots,\Sigma_\frac{N}{2}$ of the  block-circulant covariance matrix 
$
\Sigma = \text{circ}\{\Sigma_0,\Sigma_1,\dots,\Sigma_{\frac{N}{2}-1},\Sigma_{\frac{N}{2}},\Sigma_{\frac{N}{2}-1}^\top,\dots,\Sigma_1^\top\}
$
such that its inverse $\vect{S}\in\B$, i.e.
\[
   \vect{S} =\text{circ}\{S_0,S_1,\dots,S_n,0,\dots,0,S_n^\top,\dots,S_1^\top\},
\]
and the blocks $S_0,\dots,S_n$ have common support $\Omega$ as small as possible.
\end{problem}

Clearly, the matrix $\Sigma$ solving problem \ref{prob2} is the covariance of the 
reciprocal process $\vect{y}$ approximating $\vect{x}$ and featuring a sparse graphical model. 

Since we are going to identify a model for a reciprocal process, we can exploit the maximum entropy dual problem \eqref{op:mepdu} recalled before. It is worth noting that the support $\Omega$ is not known in advance, thus it has to be estimated from the data. In order to do that, inspired by \cite{SongVan}, we consider the following regularizer
\[
   h_\infty(\vect{S}) = \sum_{k>h}\,\max\left\{ |(S_0)_{hk}|,2\max_{j=1,\dots,n}|(S_j)_{hk}|,2\max_{j=1,\dots,n}|(S_j)_{kh}|\right\},
\]
which is basically a generalization of the $\ell^\infty$-norm used to induce sparsity on vectors. The optimization problem for the estimation of a sparse reciprocal model for the process is a regularized version of problem \eqref{op:mepdu}:
\begin{equation}\label{op:regop}
\begin{aligned} 
\argmin_{\vect{S}\in\B} &\quad-\log\det\vect{S} + \left<\hat{\vect{\Sigma}},\vect{S}\right>_\C + \lambda_S\,h_\infty(\vect{S})\\
\text{subject to } &\quad \vect{S}>0,
\end{aligned}
\end{equation}
where $\lambda_S>0$ is the regularization parameter. Further research is needed to understand weather \eqref{op:regop} can be seen as the dual of some kind of entropy-related optimization problem.

Notice that, although the objective function in \eqref{op:regop} is strictly convex in $\vect{S}$, it is non-differentiable due to the presence of the regularizer $h_\infty$. For this reason, we consider the dual of problem \eqref{op:regop}: as we will see below, the dual objective function is smooth and therefore it is suitable to be minimized by a projected gradient approach making the implementation of the optimization algorithm easy. Introducing the auxiliary variable $\vect{Y}\in\B$, problem \eqref{op:regop} can be rewritten as 
\begin{equation}\label{op:regop1}
\begin{aligned} 
\argmin_{\vect{S}\in\C,\vect{Y}\in\B} &\quad-\log\det\vect{S} + \left<\hat{\vect{\Sigma}},\vect{S}\right>_\C + \lambda_S\,h_\infty(\vect{Y})\\
\text{subject to } &\quad \vect{S}>0,\quad \vect{Y} = \vect{S}.
\end{aligned}
\end{equation}
Exploiting strong-duality between \eqref{op:regop1} and its dual, we address problem \eqref{op:regop1} using Lagrange multipliers theory. The Lagrangian for this problem is
\begin{align*}
\mathcal{L}(\vect{S},\vect{Y},\vect{Z}) = &-\log\det\vect{S} + \left<\hat{\vect{\Sigma}},\,\vect{S}\right>_\C + \lambda_S\,h_\infty(\vect{Y})+\left< \vect{Z},\vect{S}-\vect{Y}\right>_\C\\
= &-\log\det\vect{S} + \left<\hat{\vect{\Sigma}}+\vect{Z},\,\vect{S}\right>_\C + \lambda_S\,h_\infty(\vect{Y}) - \left<\vect{Z},\,\vect{Y}\right>_\C
\end{align*}
where $\vect{Z}\in\C$ is the Lagrange multiplier. The dual objective function is the infimum over $\vect{S}$ and $\vect{Y}$ of the Lagrangian. The unique term on $\mathcal{L}$ that depends on $\vect{Y}$ is
$ \lambda_S\,h_\infty(\vect{Y}) - \left<\vect{Z},\,\vect{Y}\right>_\C$. The latter
is bounded below if and only if
\begin{align}
	&\text{diag}(Z_j)=0,\quad  j=0,\,\dots,\,n, \label{eq:minY1}\\
	&2|(Z_0)_{kh}|+\sum_{j=1}^n\,|(Z_j)_{kh}|+|(Z_j)_{hk}|\le\frac{\lambda_S}{N},\quad k>h,\label{eq:minY2}
\end{align}
in which case the infimum is zero. Accordingly,
\[
\inf_{\vect{Y}\in\B}\,\mathcal{L}=
\left\{
\begin{split}
&-\log\det\vect{S} + \left<\hat{\vect{\Sigma}}+\vect{Z},\,\vect{S}\right>_\C \quad\text{ if }\quad\eqref{eq:minY1},\,\eqref{eq:minY2}\text{ hold,}\\\\
&-\infty\quad\text{ otherwise.}
\end{split}
\right.
\]
If \eqref{eq:minY1} and \eqref{eq:minY2} hold, it remains to minimize the strictly convex function (of $\vect{S}$)
$
  \bar{\mathcal{L}}(\vect{S}) := -\log\det\vect{S} + \left<\hat{\vect{\Sigma}}+\vect{Z},\,\vect{S}\right>_\C
$
over the symmetric, positive definite, banded block-circulant matrices. Observe that, $\forall\,\vect{Z}\in\C$, and for any sequence $\vect{S}_k>0$ converging to a singular matrix, 
\[
\lim_{k\rightarrow \infty}\bar{\mathcal{L}}(\vect{S}_k)=\infty.
\]
Accordingly, we can assume that the solution lies in the interior of the cone so that a necessary and sufficient condition for $\vect{S}_o$ to be a minimum point for $\bar{\mathcal{L}}$ is that its first Gateaux derivative computed at $\vect{S}=\vect{S}_o$ is equal to zero in every direction $\vect{\delta}\vect{S}$, namely
\begin{equation}\label{eq:ddcond}
   \vect{\delta}\bar{\mathcal{L}}(\vect{S}_o;\vect{\delta}\vect{S}) = \tr\left[ \left(-\vect{S}_o^{-1} + \hat{\vect{\Sigma}} + \vect{Z}\right)\vect{\delta}\vect{S}\right] = 0,
   \qquad \forall\,\vect{\delta}\vect{S}\in\C.
\end{equation}
Assuming that $\vect{Z}\in\C$ is such that 
\begin{equation}\label{eq:pdcstr}
	\hat{\vect{\Sigma}}+\vect{Z} > 0,
\end{equation}
condition \eqref{eq:ddcond} is satisfied if and only if $\vect{S}_o = (\hat{\vect{\Sigma}}+\vect{Z})^{-1}$. Finally, we have that
\[
\inf_{\vect{Y}\in\B,\vect{S}\in\C}\,\mathcal{L}=
\left\{
\begin{split}
&\log\det(\hat{\vect{\Sigma}}+\vect{Z}) + mN,  && \text{ if } \eqref{eq:minY1},\,\eqref{eq:minY2},\,\eqref{eq:pdcstr}\text{ hold,}\\\\
&-\infty&&\text{ otherwise.}
\end{split}
\right.
\]
The dual problem of problem \eqref{op:regop} follows straightforward
\begin{equation}\label{op:dual}
\begin{aligned} 
\argmin_{\vect{Z}\in\C} &\quad -\log\det(\hat{\vect{\Sigma}}+\vect{Z}) - mN\\
\text{subject to } &\quad \eqref{eq:minY1},\,\eqref{eq:minY2}\\
&\quad \hat{\vect{\Sigma}}+\vect{Z}>0.
\end{aligned}
\end{equation}

\begin{proposition}
	Under the assumption that $\hat{\vect{\Sigma}}\in\B$ and $\hat{\vect{\Sigma}}>0$, problem \eqref{op:dual} admits a unique solution.\\
	\begin{IEEEproof}
		Define $f(\vect{Z}):=\log\det(\hat{\vect{\Sigma}}+\vect{Z})$. Let
		\[
			\Q :=\left\{ \vect{Z}\in\C\,|\,\hat{\vect{\Sigma}}+\vect{Z}>0,\text{ and }\eqref{eq:minY1},\,\eqref{eq:minY2}\text{ hold}\right\}
		\]
		be the set of constraints of problem \eqref{op:dual}. First of all, notice that constraints \eqref{eq:minY1} and \eqref{eq:minY2} ensure that $\Q$ is a bounded subset of $\C$. Indeed, the entries of any $\vect{Z}\in\Q$ are bounded by $\lambda_S/N$ in the element-wise max-norm of the matrix. By the equivalence of norms in finite-dimensional spaces, this implies in particular that $\|\vect{Z}\|_\C<\infty$ for any $\vect{Z}\in\Q$. Let now $(\vect{Z}^{(k)})_{k\in\NN}$ be a generic sequence of elements of $\Q$ converging to some $\bar{\vect{Z}}\in\C$, such that $\hat{\vect{\Sigma}}+\bar{\vect{Z}}\ge0$ singular. Then
		\[
			\lim_{k\to\infty}\,-\log\det(\hat{\vect{\Sigma}}+\vect{Z}^{(k)})=+\infty,
		\]
		and therefore $\vect{Z}^{(k)}$ is not an infimizing sequence. Hence, we can restrict our attention to the closed subset of $\Q$ defined by 
		\[
		\bar{\Q}:=\left\{ \vect{Z}\in\C\,|\,\hat{\vect{\Sigma}}+\vect{Z}\ge\epsilon I_{mN},\text{ and }\eqref{eq:minY1},\,\eqref{eq:minY2}\text{ hold}\right\}
		\]
		with $\epsilon>0$ small enough. By what we have shown till now, the function $f$ is continuous on the compact set $\bar{\Q}$ and therefore admits at least one minimum point. Since $f$ is strictly convex, the minimum is unique.
	\end{IEEEproof}
\end{proposition}

\begin{proposition}
	Under the assumption that $\hat{\vect{\Sigma}}\in\B$ and $\hat{\vect{\Sigma}}>0$, problem \eqref{op:regop} admits a unique solution $\vect{S}_o$.
\end{proposition}
\begin{IEEEproof}
	Notice that problem \eqref{op:regop} is a strictly feasible (for instance, pick $X=I_{mN}$) convex optimization problem. Accordingly, Slater's condition holds, hence strong duality holds between \eqref{op:regop} and its dual. The strong duality between problems \eqref{op:regop} and \eqref{op:dual} and the existence of a unique optimum $\vect{Z}_o$ for the dual problem \eqref{op:dual}, imply that there exists a unique $\vect{S}_o\in\C$ so that
	$
	\vect{S}_o = \left(\hat{\vect{\Sigma}} + \vect{Z}_o\right)^{-1}
	$
	which solves the primal problem \eqref{op:regop}.
\end{IEEEproof}

\section{NUMERICAL EXAMPLE}\label{sec:ne}
\begin{figure}[h!]
		\begin{minipage}[.c]{0.15\textwidth}
			\centering
			\includegraphics[scale=0.24]{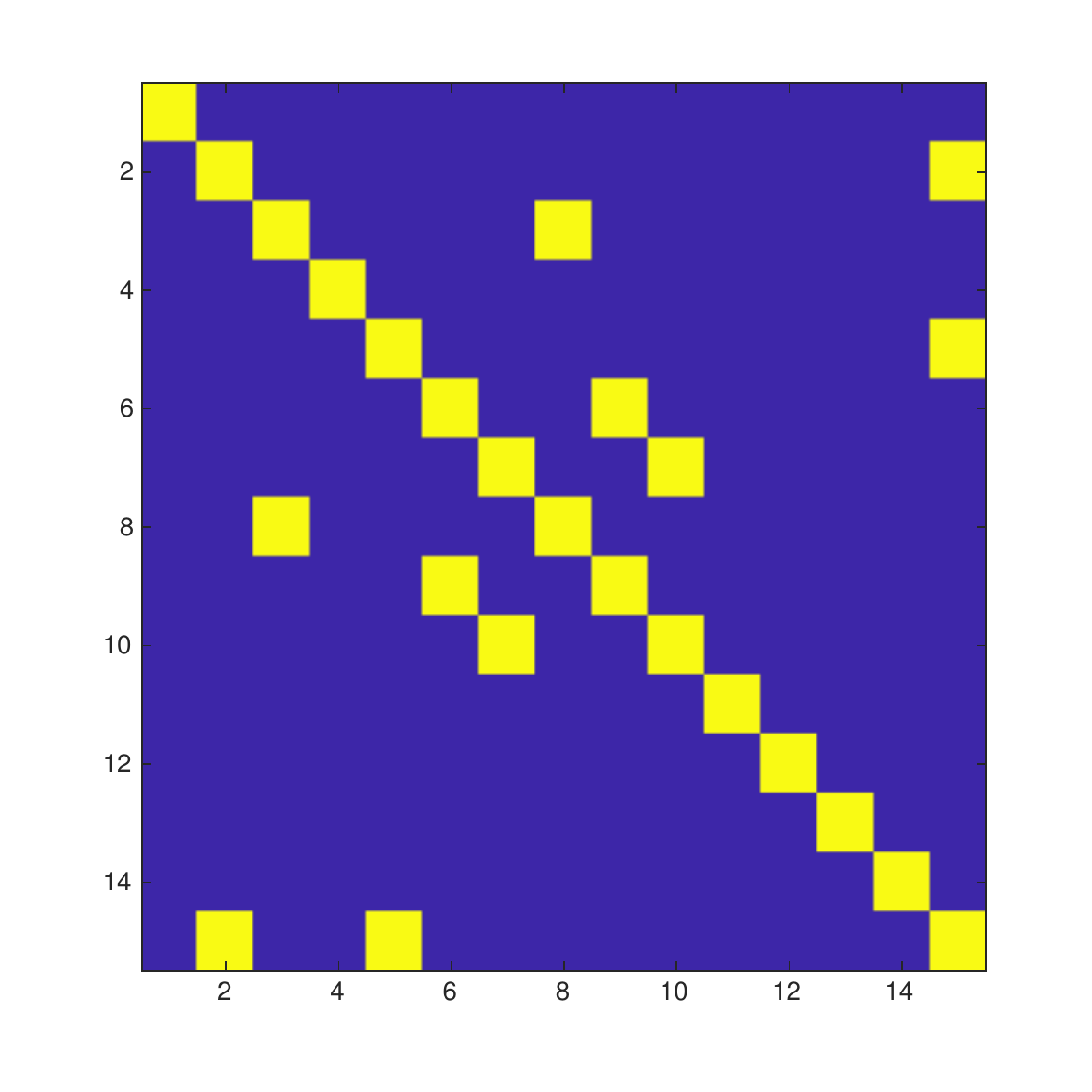}
		\end{minipage}%
		\begin{minipage}[.c]{0.15\textwidth}
			\centering
			\includegraphics[scale=0.24]{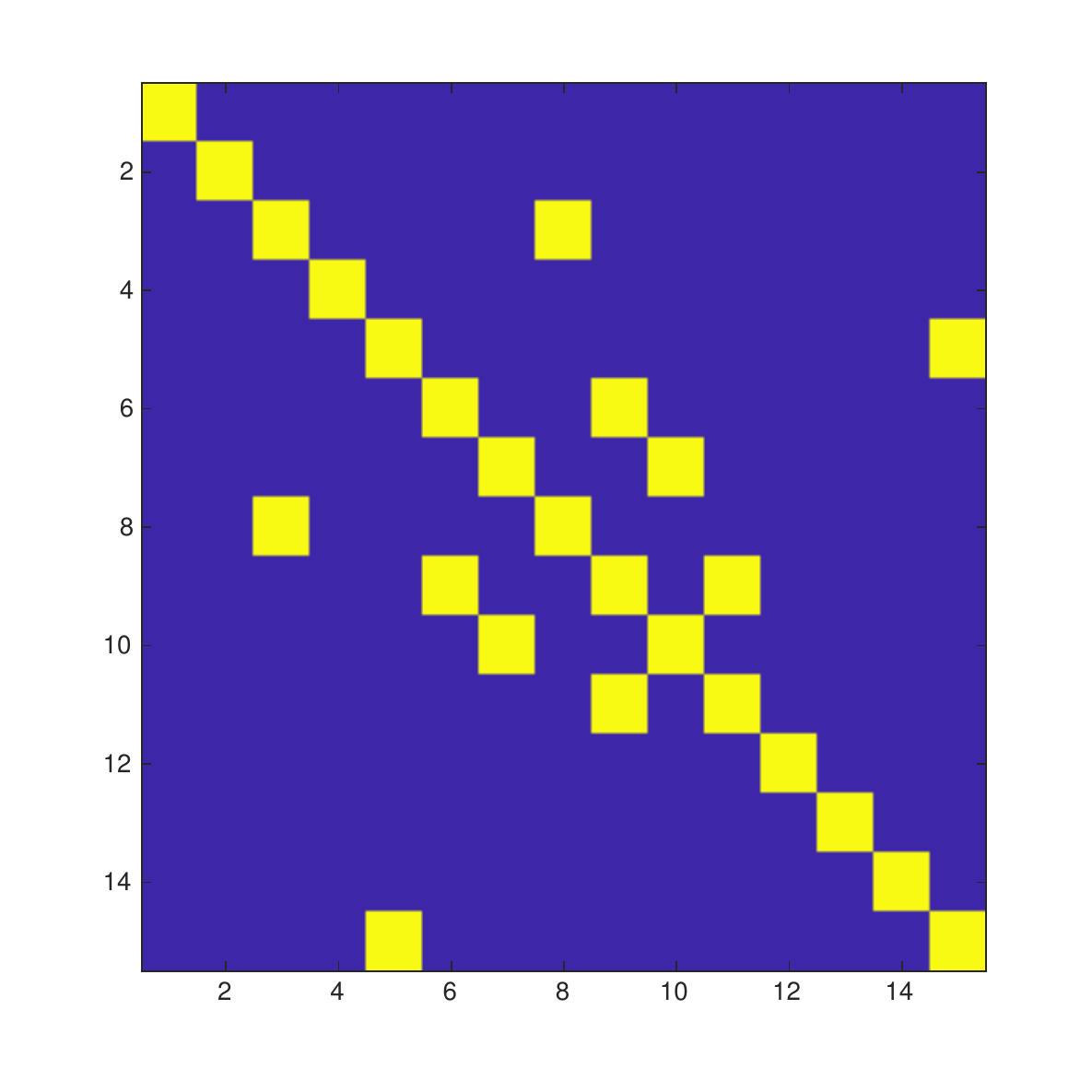}
		\end{minipage}%
		\begin{minipage}[.c]{0.15\textwidth}
			\centering
			\includegraphics[width=3.0cm,height=2.9cm]{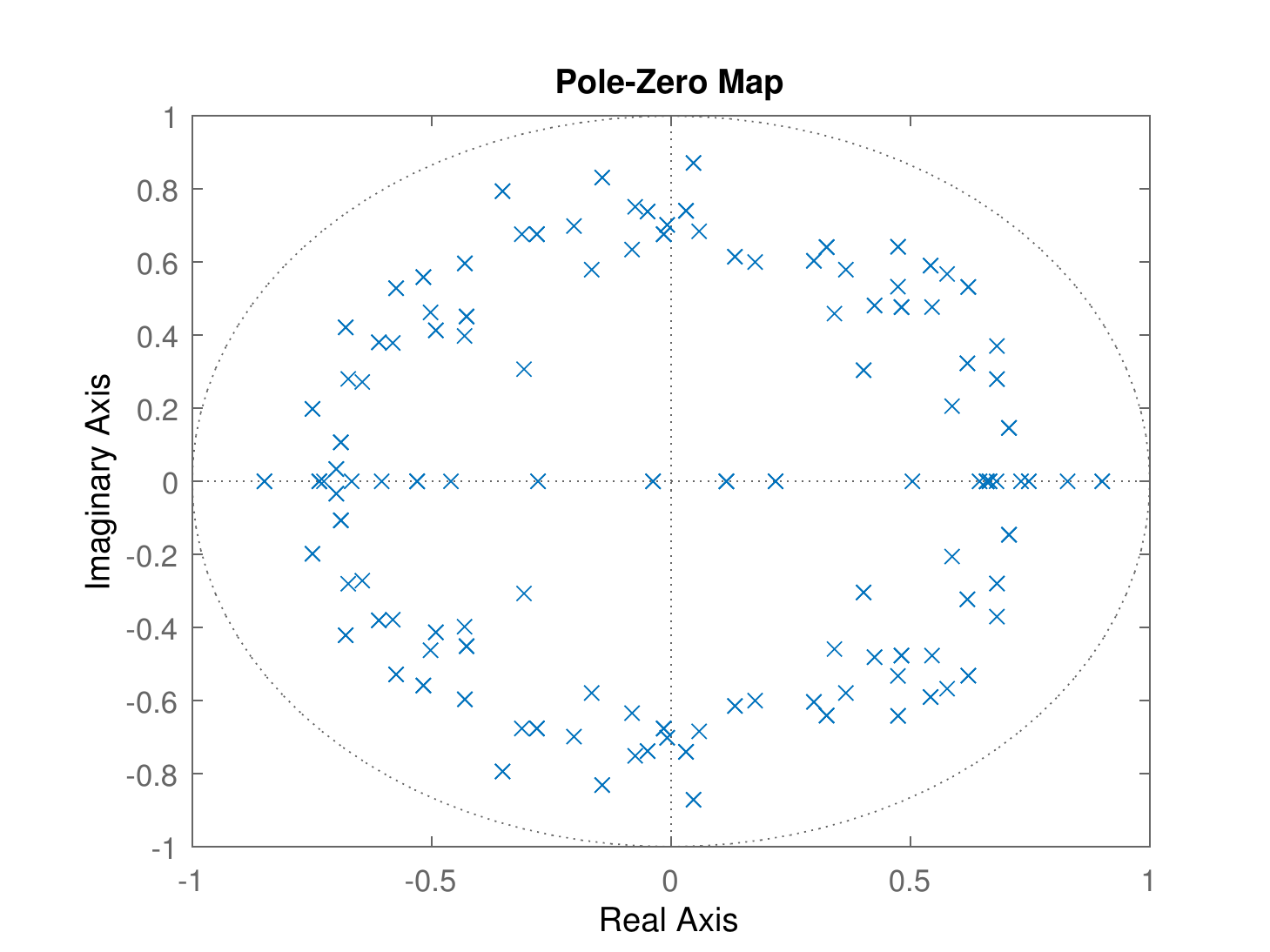}
		\end{minipage}
		\caption{Sparsity pattern of the inverse spectrum: true (left) and estimated (center). The yellow squares represent the conditional dependent pairs while the blue ones stand for the conditional independent pairs. Poles of the model used to generate the data (right).}
		\label{fig:omega}
\end{figure}

Here we present a numerical example that shows how the algorithm works in practice. We consider the AR model $\vect{x}(t)=A(z)\vect{x}(t)+\vect{w}(t)$ with $A(z)=\sum_{k=1}^n A_k z^{-k}$, $m=15$, $n=8$, $\vect{w}(t)$ is white Gaussian noise with variance equal to $15.9$ and the poles $p_j$ of the shaping filter $[I-A(z)]^{-1}$ are depicted in Figure \ref{fig:omega} (right) and they are such that $|p_j|\le 0.9$. Figure \ref{fig:omega} also shows the sparsity pattern of the true inverse spectrum with the sparsity pattern estimated by the proposed algorithm when the regularization parameter is $\lambda_S=125$ and $T=1000$ samples are used to estimate the covariance lags $\hat{R}_k$.
\\
From Figure \ref{fig:omega} we can see that only one non-zero element has been estimated in a wrong position. The average relative error between the estimated and the true spectra is about $4\%$.

\section{CONCLUSIONS AND FUTURE DEVELOPMENTS}\label{sec:concl}
In this paper we presented an identification procedure  for a sparse graphical model associated with a reciprocal process. 
As discussed in detail in \cite{PICCI_LINDQUIST} and recalled in the introduction, the reciprocal approximation must be understood as an attempt in the direction of the development of an effective procedure for the identification of ARMA graphical models that scales robustly with the product of the process dimension by the length of an accurate AR approximation of the original process. This is a promising theoretical idea that will be tested in simulation and in real examples. Moreover, this approach can be pushed forward in many directions: for example there is the possibility of adding a (small) number of latent variables 
to the picture in order  to provide a better approximation of the dynamics of the original process.



\begin{thebibliography}{10}

\bibitem{BDS}
M.~S. Chen~M. and L.~Y., ``Big data: A survey,'' {\em Mobile Netw Appl},
  no.~19, pp.~171--209, 2014.

\bibitem{Lauritzen}
S.~Lauritzen, {\em Graphical Models}.
\newblock Oxford, U.K.: Oxford university press, 1996.

\bibitem{candes2012exact}
E.~Candes and B.~Recht, ``Exact matrix completion via convex optimization,''
  {\em Communications of the ACM}, vol.~55, no.~6, pp.~111--119, 2012.

\bibitem{candes2010matrix}
E.~Candes and Y.~Plan, ``Matrix completion with noise,'' {\em Proceedings of
  the IEEE}, vol.~98, no.~6, pp.~925--936, 2010.

\bibitem{ZorzCh}
M.~Zorzi and A.~Chiuso, ``Sparse plus low rank network identification: A
  nonparametric approach,'' {\em Automatica}, vol.~76, no.~2, pp.~355--366,
  2017.

\bibitem{songsiri2010graphical}
J.~Songsiri, J.~Dahl, and L.~Vandenberghe, ``Graphical models of autoregressive
  processes,'' {\em Convex optimization in signal processing and
  communications}, pp.~89--116, 2010.

\bibitem{SongVan}
J.~Songsiri and L.~Vandenberghe, ``Topology selection in graphical models of
  autoregressive processes,'' {\em J. Mach. Learn. Res.}, vol.~11,
  pp.~2671--2705, 2010.

\bibitem{LindqARMA}
E.~Avventi, A.~Lindquist, and B.~Wahlberg, ``Arma identification of graphical
  models,'' {\em IEEE Transactions on Automatic Control}, vol.~58,
  pp.~1167--1178, May 2013.

\bibitem{PICCI_LINDQUIST}
A.~Lindquist and G.~Picci, ``The circulant rational covariance extension
  problem: The complete solution,'' {\em IEEE Transactions on Automatic
  Control}, vol.~58, pp.~2848--2861, Nov 2013.

\bibitem{BEL01}
C.~I. Byrnes, P.~Enqvist, and A.~Lindquist, ``Cepstral coefficients, covariance
  lags, and pole-zero models for finite data strings,'' {\em IEEE Transactions
  on Signal Processing}, vol.~49, pp.~677--693, Apr 2001.

\bibitem{BEL02}
C.~I. Byrnes, P.~Enqvist, and A.~Lindquist, ``Identifiability and
  well-posedness of shaping filter parametrizations: A global analysis
  approach,'' {\em SIAM Journal on Control and Optimization}, vol.~41, no.~1,
  pp.~23--59, 2002.

\bibitem{MK85}
B.~R. Musicus and A.~M. Kabel, ``Maximum entropy pole-zero estimation,'' Tech.
  Rep. 510, Massachusetts Institute of Technology, Aug 1985.

\bibitem{ChaParWill}
V.~Chandrasekaran, P.~A. Parrilo, and A.~S. Willsky, ``Latent variable
  graphical model selection via convex optimization,'' {\em Ann. Statist.},
  vol.~40, pp.~1935--1967, 08 2012.

\bibitem{ZorzSep}
M.~Zorzi and R.~Sepulchre, ``{AR} identification of latent-variable graphical
  models,'' {\em IEEE Transactions on Automatic Control}, vol.~61, no.~9,
  pp.~2327--2340, 2016.

\bibitem{CarFerrPav}
F.~P. Carli, A.~Ferrante, M.~Pavon, and G.~Picci, ``A maximum entropy solution
  of the covariance extension problem for reciprocal processes,'' {\em IEEE
  Trans. on Automatic Control}, vol.~56, pp.~1999--2012, Sept 2011.

\bibitem{Rec1}
B.~C. Levy and A.~Ferrante, ``Characterization of stationary discrete-time
  gaussian reciprocal processes over a finite interval,'' {\em SIAM J. on
  Matrix Analysis and Applications}, vol.~24, no.~2, pp.~334--355, 2002.

\bibitem{Rec2}
B.~C. Levy, ``Regular and reciprocal multivariate stationary gaussian
  reciprocal processes over z are necessarily markov,'' {\em J. Math. Syst.
  Est. Control}, vol.~2, pp.~134--154, 1992.

\bibitem{Rec3}
B.~C. Levy, R.~Frezza, and A.~J. Krener, ``Modeling and estimation of
  discrete-time gaussian reciprocal processes,'' {\em IEEE Transactions on
  Automatic Control}, vol.~35, pp.~1013--1023, Sep 1990.

\bibitem{ringh2016multidimensional}
A.~Ringh, J.~Karlsson, and A.~Lindquist, ``Multidimensional rational covariance
  extension with applications to spectral estimation and image compression,''
  {\em SIAM Journal on Control and Optimization}, vol.~54, no.~4,
  pp.~1950--1982, 2016.

\bibitem{FerrLinAlg}
F.~Carli, A.~Ferrante, M.~Pavon, and G.~Picci, ``An efficient algorithm for
  maximum entropy extension of block-circulant covariance matrices,'' {\em
  Linear Algebra and its Applications}, vol.~439, no.~8, pp.~2309 -- 2329,
  2013.

\bibitem{ringh2015fast}
A.~Ringh and J.~Karlsson, ``A fast solver for the circulant rational covariance
  extension problem,'' in {\em Control Conference (ECC), 2015 European},
  pp.~727--733, IEEE, 2015.

\bibitem{tesi}
D.~Alpago, ``On the identification of sparse plus low-rank graphical models,''
  Master's thesis, University of Padova, Padova, Italy, 2017.

\bibitem{structcov1}
A.~Ferrante, M.~Pavon, and M.~Zorzi, ``A maximum entropy enhancement for a
  family of high-resolution spectral estimators,'' {\em IEEE Transactions on
  Automatic Control}, vol.~57, pp.~318--329, Feb 2012.

\bibitem{structcov2}
M.~Zorzi and A.~Ferrante, ``On the estimation of structured covariance
  matrices,'' {\em Automatica}, vol.~48, no.~9, pp.~2145 -- 2151, 2012.

\bibitem{Gray}
R.~M. Gray, ``Toeplitz and circulant matrices: A review,'' {\em Foundations and
  Trends{\textregistered} in Communications and Information Theory}, vol.~2,
  no.~3, pp.~155--239, 2006.

\bibitem{BurgPhd}
J.~Burg, {\em Maximum entropy spectral analysis}.
\newblock PhD thesis, Stanford University, Dept. of Geophysics, Stanford, CA,
  1975.

\bibitem{family_zorzi}
M.~Zorzi, ``A new family of high-resolution multivariate spectral estimators,''
  {\em IEEE Transactions on Automatic Control}, vol.~59, pp.~892--904, April
  2014.

\bibitem{ZORZI201587}
M.~Zorzi, ``An interpretation of the dual problem of the {THREE}-like
  approaches,'' {\em Automatica}, vol.~62, pp.~87 -- 92, 2015.

\end{thebibliography}
\end{document}